\documentclass[12pt,a4paper]{article}

\usepackage{amsmath,amssymb}
\usepackage{txfonts}
\usepackage{bm}

\advance \textwidth -60truemm\relax

\advance\oddsidemargin -1truein\relax

\evensidemargin=\oddsidemargin

\advance\textheight -60truemm\relax
\advance\topmargin -1truein\relax
\unitlength\textwidth
\divide\unitlength by 150\relax

\def\RRR{\mathbb{R}}
\def\CCC{\mathbb{C}}
\def\FFF{\mathbb{F}}
\def\KKK{\mathbb{K}}
\bmdefine{\aaa}{a}
\bmdefine{\bbb}{b}
\bmdefine{\uuu}{u}
\bmdefine{\vvv}{v}
\bmdefine{\www}{w}
\bmdefine{\eee}{e}
\bmdefine{\xxx}{x}
\bmdefine{\yyy}{y}
\bmdefine{\zerovec}{0}

\newcommand{\rank}{\mathrm{rank}}
\newcommand{\diag}{{\mathrm{Diag}}}

\newcommand{\maxrank}{{\mathrm{max.rank}}}

\newcommand{\brk}{{\mathrm{brk}}}
\newcommand{\maxbrk}{{\mathrm{max.brk}}}

\newcommand{\jaja}{JaJa}
\newcommand{\card}{\mathrm{Card}}

\numberwithin{equation}{section}
\newtheorem{thm}[equation]{Theorem}
\newtheorem{lemma}[equation]{Lemma}
\newtheorem{cor}[equation]{Corollary}
\newtheorem{prop}[equation]{Proposition}
\newtheorem{example}[equation]{Example}
\newtheorem{definition}[equation]{Definition}

\makeatletter
\def\proof{\@ifnextchar[{\proof@}{\proof@[]}}
\def\proof@[#1]{\noindent{\bf Proof #1}\ }

\makeatother

\begin{document}

\title{Rank of $3$-tensors with $2$ slices and Kronecker canonical forms}
\author{Toshio Sumi, Mitsuhiro Miyazaki and Toshio Sakata%
\thanks{Kyushu University, Kyoto University of Education and Kyushu University}
\thanks{The authors were supported partially by Grant-in-Aid for Scientific Research (B) (No.~20340021) of the Japan Society for the Promotion of Science.}
}
\maketitle
\pagestyle{plain}


\begin{abstract}
Tensor type data are becoming important recently in various application 
fields.  
We determine a rank of a tensor $T$ so that
$A+T$ is diagonalizable for a given $3$-tensor $A$
with $2$ slices over the complex and real number field.
\end{abstract}

%
\section{Introduction}
Tensor type data are becoming important recently in various application 
fields (for example see Miwakeichi et al.~\cite{Miwa:2004}, Vasilescu and Terzopoulos~\cite{Vasilescu-Terzopoulos:2005} and Muti and Bourennane~\cite{Muti-Bourennane:2007}).   
The factorization of a tensor to a sum of rank $1$ tensors means that the data is expressed by a sum of data with simplest structure, and 
we may have better understanding of data. 
This is an essential attitude for data analysis and therefore 
the problem of tensor factorization is an essential one for applications.  
In this paper we consider the rank problem of $3$-tensors with $2$ slices. 
This was studied in the 1970's and 1980's by many authors.
\jaja{} \cite{JaJa:1979b} 
gave the rank for a $3$-tensors with $2$ slices.  He used Kronecker canonical forms of the pencil of two matrices.
Results by Brockett and Dobkin \cite{Brockett-Dobkin:1973, Brockett-Dobkin:1978} are useful for giving a lower bound.
\jaja{} showed that the rank of a Kronecker canonical form without
regular pencils is equal to the sum of the ranks of direct summand.
However, the rank of a Kronecker canonical form is not equal to
the sum of the ranks of direct summand in general and it depends
on invariant polynomials.
This causes to be difficult to determine the rank of tensors.
Our aim is to determine a rank of a tensor $T$ so that
$A+T$ is diagonalizable for a given $3$-tensor $A$
with $2$ slices (see Theorem~\ref{thm:diagonal}). 
In this paper we consider ranks of tensors over the complex and 
real number field.

\section{Kronecker canonical forms}

We consider the case of the complex number field and the real number field. 
Let $\mathbb{F}=\mathbb{R},\mathbb{C}$. We note that 
$(A_1,\ldots,A_r)$ denotes a horizontally posed $m \times nr$ matrix 
and $(A_1;\ldots;A_r)$ denotes a tensor whose $k$-th slice is an 
$m\times n$ matrix  $A_k$.
A tensor $(A;B)$ is called {\sl diagonalizable} if
there are 
an $m\times m$ nonsingular matrix $P$, an $n\times n$ nonsingular 
matrix $Q$ and diagonal matrices $D_A$, $D_B$ such that
$PAQ=(D_A,O)$, $PBQ=(D_B,O)$ for $m\leq n$ 
and $PAQ=(D_A,O)^T$, $PBQ=(D_B,O)^T$ for $m>n$.
%
Note that if $A$ is an $n\times n$ matrix, then $(E_n;A)$ is
diagonalizable if and only if $A$ is similar to a diagonal matrix,
i.e., 
there is a nonsingular matrix $P$ such that $PAP^{-1}$ is a 
diagonal matrix.

We summarize briefly about Kronecker canonical forms.

\begin{lemma}[{\cite[(30) in \S4, XII]{{Gantmacher:1959}}}] \label{lem:KroneckerCanonicalForm}
Let $A$ and $B$ be  $m\times n$ 
rectangular matrices. 
Then there are nonsingular matrices $P$ and $Q$ such that 
$$P(A;B)Q=(PAQ;PBQ)$$ 
is of a block diagonal form
$$\diag((S_1;T_1),\ldots,(S_r;T_r)),$$
where each $(S_j;T_j)$ is one of the following  
\begin{enumerate}
\item[{\rm (A)}] 
	$k\times \ell\times 2$ tensor $(O;O)$, 
\item[{\rm (B)}] 
	$k\times k\times 2$ tensor $(\alpha E_k+J_k;E_k)$, 
\item[{\rm (C)}] 
	$2k\times 2k\times 2$ tensor $(C_k(c,s)+J_k\otimes E_2;E_{2k})$, 
	$s\ne 0$, 
\item[{\rm (D)}] 
	$k\times k\times 2$ tensor $(E_k;J_k)$, 
\item[{\rm (E)}] 
	$k\times (k+1)\times 2$ tensor $((\zerovec,E_k);(E_k,\zerovec))$, 
\item[{\rm (F)}] 
	$(k+1)\times k\times 2$ tensor 
	$(\begin{pmatrix} \zerovec^T\cr E_k\end{pmatrix};
	  \begin{pmatrix} E_k\cr \zerovec^T\end{pmatrix})$. 
\end{enumerate}
Here $E_k$ is the $k\times k$ identity matrix, 
$J_k=\begin{pmatrix} 0&1&&O \cr 
  \vdots&\ddots&\ddots \cr \vdots&&\ddots&1\cr
  0&\cdots&\cdots&0 \end{pmatrix}$
is a $k \times k $ square matrix and
$C_k(c,s)=E_k\otimes \begin{pmatrix} c&-s\cr s&c\end{pmatrix}=
  \diag(\begin{pmatrix} c&-s\cr s&c\end{pmatrix},\ldots ,
  \begin{pmatrix} c&-s\cr s&c\end{pmatrix})$ 
is a $2k \times 2k$ square matrix. 
\end{lemma}

This decomposition is called the Kronecker canonical form.
It is unique up to permutations of blocks.
Note that tensors of type \eqref{typeA} include ones when $k>0$ and $\ell=0$, 
or $k=0$ and $\ell>0$, where a direct sum of a $0\times \ell$ tensor of 
type \eqref{typeA}
and an $s\times t$ tensor $(X;Y)$ means a $k\times (\ell+t)$ tensor $((O,X),(O,Y))$.
Also note that type \eqref{typeC} does not appear over the complex number field $\mathbb{C}$, and if $\alpha$ is not real in \eqref{typeB}, then type \eqref{typeC} appears over the real number field $\mathbb{R}$.
\par

First we note the following 

\begin{lemma}\label{lem:diag}
Let $A$ be an $\ell\times \ell$ matrix with entries in $\FFF$.
Then 
$\rank_\FFF(E_\ell;A)=\ell$ if and only if 
$(E_\ell;A)$
is diagonalizable over $\FFF$.
\end{lemma}

To estimate the ranks of tensors of types \eqref{typeB}, \eqref{typeC} and \eqref{typeD},
we recall some basic facts of linear algebra.

\begin{definition}\label{def:compmat}\rm
Let $f(x)=x^n+a_1x^{n-1}+\cdots+a_{n-1}x+a_n$  be a monic polynomial
with coefficients in $\FFF$.
The matrix
$$
M:=\begin{pmatrix}
&&& -a_n \\
1&&& -a_{n-1} \\
&\ddots &&\vdots \\
&&1&-a_1
\end{pmatrix}
$$
is called a {\sl companion matrix} for $f(x)$.
\end{definition}
Note that $f(x)$ is both the characteristic polynomial and the minimal polynomial of $M$.
For any monic polynomial $g(x)=x^n+b_1x^{n-1}+\cdots+b_{n-1}x+b_n$
of degree $n$, set
$$
N:=\begin{pmatrix}
0&\cdots&0& b_n-a_n \\
0&\cdots&0&b_{n-1} -a_{n-1} \\
\vdots& &\vdots&\vdots \\
0&\cdots&0&b_1-a_1
\end{pmatrix} \enspace.
$$
Then $M-N$ is the companion matrix for $g(x)$.
In particular, by taking $g(x)$ to be 
a product of distinct monic linear polynomials,
we see the following
\begin{lemma}\label{lem:+1}
For any companion matrix $M$ with entries in $\FFF$,
there is a tensor $(O;N)$ with rank at most $1$ 
such that
$(E_n;M)-(O;N)$ is diagonalizable.
In particular, it holds
$\rank_{\FFF}(E_n;M)\leq n+1$.
\end{lemma}

The following three lemmas are well known and easily proved.

\begin{lemma}\label{lem:blockJordan}
Let $A$ be an $\ell\times \ell$ matrix.
If the minimal polynomial $h(x)$ of $A$ has no 
multiple factor, then the minimal polynomial of
$E_k\otimes A+J_k\otimes E_\ell$ is $h(x)^k$.
\end{lemma}

\begin{lemma}\label{lem:minofdiag}
The minimal polynomial of $\diag(A,B)$ is the
least common multiple of the minimal polynomials of $A$
and $B$.
\end{lemma}

\begin{lemma}\label{lem:min=char}
Let $A$ be a square matrix whose minimal polynomial
is equal to the characteristic polynomial.
Then $A$ is similar to the companion matrix for the
characteristic polynomial of $A$, 
i.e.,
there is a nonsingular polynomial $P$ such that
$PAP^{-1}$ is the companion matrix for the characteristic
polynomial of $A$.
\end{lemma}

By Lemma \ref{lem:+1}, we see the following

\begin{cor}\label{cor:+1}
Suppose an $\ell\times \ell$ matrix $A$ satisfies the condition
of Lemma \ref{lem:min=char}.
Then there is an $\ell\times \ell\times 2$ tensor $T$ with
rank at most $1$ such that $(E_{\ell};A)-T$ is diagonalizable.
In particular,
$\rank_{\FFF}(E_\ell;A)\leq \ell+1$.
\end{cor}

Now we examine the tensors of types \eqref{typeB}, \eqref{typeC} and \eqref{typeD}.
\begin{lemma}\label{lem:BCD}
For an $\ell\times \ell\times 2$ tensor $T$ of type \eqref{typeB},
\eqref{typeC} or \eqref{typeD},
there is a tensor $T'$ with rank at most $1$ such that
$T-T'$ is diagonalizable.
In particular, 
$\rank_{\FFF}(T)\leq \ell+1$.
\end{lemma}

\begin{proof}
First consider the tensor of type \eqref{typeB}.
Since the minimal polynomial of $1\times 1$ matrix $(\alpha)$
is $x-\alpha$,
the minimal polynomial of $\alpha E_\ell+J_\ell$ is $(x-\alpha)^\ell$
by Lemma \ref{lem:blockJordan}.
So the minimal polynomial of $\alpha E_\ell+J_\ell$ is equal to the
characteristic polynomial of it.
Therefore the result follows by Corollary \ref{cor:+1}.

Type \eqref{typeD} is a special case of type \eqref{typeB}.

Finally, we consider a tensor of type \eqref{typeC}.
Note that $\FFF=\RRR$ in this case.
Since the minimal polynomial of $C_1(c,s)$ is 
an irreducible polynomial of degree 2, the result 
follows from Lemma \ref{lem:blockJordan} and Corollary \ref{cor:+1}.
\end{proof}

\jaja{} showed that 
$\rank_{\FFF}((\zerovec,E_k);(E_k,\zerovec))=k+1$ 
(see \cite[Theorem~2.1]{JaJa:1979b}).
The following is obtained from his proof.

\begin{lemma} \label{lem:GSDcaseE&F} 
For a $k\times (k+1)\times 2$ tensor 
$(A;B)=((\zerovec,E_k);(E_k,\zerovec))$,  
there are a rank $1$ matrix $M$, a nonsingular matrices $P,Q$ and 
numbers $s_1,\ldots s_k$ and such that 
$$P((A;B)-(M;O))Q=(\diag(s_1,\ldots,s_k),\zerovec);(E_k,\zerovec))
\enspace.$$ 
\end{lemma}

Note that the similar result as above holds 
for a $(k+1)\times k\times 2$ tensor 
$(A;B)=(\begin{pmatrix} \zerovec^T\cr E_k\end{pmatrix};
  \begin{pmatrix} E_k\cr \zerovec^T\end{pmatrix})$
since $(A;B)=((\zerovec,E_k);(E_k,\zerovec))^T$.

We denote by $\langle A_1,\ldots, A_m\rangle$ a vector space spanned by
matrices $A_1,\ldots,A_m$.

\begin{thm} \label{thm:x2mn>}
Let $m\leq n\leq 2m$ and $(A_1;\ldots;A_m)$ be a $2\times n\times m$ tensor.
Suppose that $\dim\langle A_1,\ldots,A_m\rangle=m$.
Let $\ell\leq\lfloor n/2\rfloor$ be an integer 
such that
$\rank({A'_1}^T, \ldots, {A'_\ell}^T)=2\ell$
for any $A'_j=A_j+c_{j,\ell+1}A_{\ell+1}+\cdots+c_{j,m}A_m$, 
$1\leq j\leq \ell$.
Then
$\rank_{\FFF}(A_1;\cdots;A_m)\geq m+\ell$.
\end{thm}

\begin{proof}
Since $\dim\langle A_1,\ldots,A_m\rangle=m$, it holds that
$\rank_{\FFF}(A_1;\ldots;A_m)\geq m$.
Assuming that each $A_j$ for $1\leq j\leq m$ is a linear combination 
of rank $1$ matrices $C_1, \ldots, C_{m+q}$,
we derive that $q$ must be larger than or equal to $\ell$.
Putting 
$$
A_i=\sum_{j=1}^{m+q}\alpha_{ij}C_j, \quad i=1, \ldots, m
$$
since the $m\times (m+q)$ matrix $(\alpha_{ij})$ has rank $m$, if necessary, 
exchanging suffixes, without loss of generality, we can assume the first $m$ columns of $(\alpha_{ij})$ are linearly independent. 
Let $(\beta_{ij})$ be its inverse matrix and take $i_1, i_2, \ldots, i_m$ so that $(\beta_{i_s,t})_{1\leq s,t\leq k}$ become nonsingular matrices for all $1\leq k\leq m$.  
Putting 
$$B_s=\sum_{j=1}^{m}\beta_{i_s,j}A_j,
\quad s=1,\ldots,\ell
$$
and define a $2\ell\times n$ matrix $X$ by 
$X=\begin{pmatrix} B_1\cr \vdots\cr B_\ell\end{pmatrix}$.
Let $P$ be the inverse matrix of the $\ell\times \ell$ square matrix $(\beta_{i_s,t})_{1\leq s,t\leq \ell}$.
Since $P(\beta_{i_s,j})=(E_\ell,\ast)$ for the $\ell\times m$ matrix $(\beta_{i_s,j})$,
we have  
$$(P\otimes E_2)X=\begin{pmatrix}
A_1+(\mbox{a linear combination of $A_{\ell+1}, \ldots, A_m$})\\
A_2+(\mbox{a linear combination of $A_{\ell+1}, \ldots, A_m$})\\
\vdots\\
A_\ell+(\mbox{a linear combination of $A_{\ell+1}, \ldots, A_m$})
\end{pmatrix}\enspace.$$
Then $\rank(X)=2\ell$ by assumption. 
On the other hand, since 
\begin{eqnarray*}
B_s&=&\sum_{j=1}^m\beta_{i_s, j}A_j 
  =\sum_{j=1}^m \beta_{i_s,j}\sum_{k=1}^{m+q}\alpha_{jk}C_k
  =\sum_{k=1}^{m+q}(\sum_{j=1}^m \beta_{i_s,j}\alpha_{jk})C_k\\
&=&C_{i_s}+\sum_{k=1}^{q} \gamma_{sk} C_{m+k} \\
\end{eqnarray*}
where $\gamma_{sk}=\displaystyle\sum_{j=1}^m \beta_{i_s,j}\alpha_{j,m+k}$ for $s=1,\ldots,\ell$ and $k=1,\ldots, q$, we have 
$$
X=\begin{pmatrix} 1 \cr 0\cr \vdots \cr 0 \end{pmatrix} \otimes C_{i_1} 
+\cdots
+\begin{pmatrix} 0 \cr \vdots\cr 0 \cr 1 \end{pmatrix} \otimes C_{i_\ell} 
+ \begin{pmatrix} \gamma_{11} \cr \gamma_{21}\cr \vdots \cr \gamma_{\ell 1} \end{pmatrix} \otimes C_{m+1} 
+\cdots
+ \begin{pmatrix} \gamma_{1q} \cr \gamma_{2q}\cr \vdots \cr \gamma_{\ell q} \end{pmatrix} \otimes C_{m+q} \\
$$
and therefore $X$ becomes a linear combination of $\ell+q$ matrices of rank $1$.  This means that $q\geq \ell$, which completes the proof.
\end{proof}

\begin{cor} \label{cor:x2mn>}
Let $m\leq n\leq 2m$ and $\ell\leq \lfloor n/2\rfloor$.
Let $X_{11}$, $X_{22}$ and $Y$ be nonsingular $(n-\ell)\times (n-\ell)$, 
$(m+\ell-n)\times (m+\ell-n)$ and $\ell\times\ell$ matrices respectively.
We define $m\times n$ matrices $A$ and $B$ by 
$$
A=\begin{pmatrix} X_{11}&X_{12}&O\\ O&X_{22}&O\end{pmatrix},\quad
B=\begin{pmatrix} O&Y\cr O&O\end{pmatrix} \enspace.
$$
Then
$\rank_{\FFF}(A;B)= m+\ell$.
\end{cor}

\begin{proof}
Set $X_{11}=\begin{pmatrix} \xxx_{11}^T\\ \vdots\\ \xxx_{n-\ell,1}^T\end{pmatrix}$,
$X_{12}=\begin{pmatrix} \xxx_{12}^T\\ \vdots\\ \xxx_{n-\ell,2}^T\end{pmatrix}$,
$X_{22}=\begin{pmatrix} \xxx_{n-\ell+1,2}^T\\ \vdots\\ \xxx_{m,2}^T\end{pmatrix}$.
Let $Y=(Y_1,Y_2)$, where 
$Y_1=\begin{pmatrix} \yyy_{11}^T\\ \vdots\\ \yyy_{\ell,1}^T\end{pmatrix}$ 
is a $\ell\times (m+\ell-n)$ matrix and
$Y_2=\begin{pmatrix} \yyy_{12}^T\\ \vdots\\ \yyy_{\ell,2}^T\end{pmatrix}$
is a $\ell\times(n-m)$ matrix.
We take $(A;B)$ as an array with $m$ slices of $2 \times n$ matrices 
$A_1, A_2, \ldots, A_m$:
$$A_i=\begin{pmatrix}
\xxx_{i1}^T&\xxx_{i2}^T&\zerovec^T\\
\zerovec^T&\yyy_{i1}^T&\yyy_{i2}^T
\end{pmatrix}
$$
for $1\leq i\leq n-\ell$ and
$$A_i=\begin{pmatrix}
\zerovec^T&\xxx_{i2}^T&\zerovec^T\\
\zerovec^T&\zerovec^T&\zerovec^T
\end{pmatrix}
$$
for $n-\ell+1\leq i\leq m$.
Here $\yyy_{i1}=\yyy_{i2}=\zerovec$ if $i>\ell$.
Since $\rank(A)=m$, it holds $\dim\langle A_1,\ldots,A_m\rangle=m$
and also by assumption $A_1,\ldots,A_m$ satisfy the assumption of 
Theorem~\ref{thm:x2mn>} and then $\rank_{\FFF}(A;B)\geq m+\ell$.
Conversely, we have $\rank_{\FFF}(A;B)\leq \rank(A)+\rank(B)=m+\ell$.
\end{proof}

\begin{example} \label{ex:maxrank} For the tensor 
$X=((E_m,O);\begin{pmatrix} O&E_{\lfloor n/2\rfloor}\cr O&O\end{pmatrix})$ 
of $\FFF^{m\times n\times 2}$ with $m\leq n\leq 2m$, 
it holds that
$\rank_{\FFF}(X)=m+\lfloor n/2\rfloor$.
\end{example}

\begin{thm} \label{thm:JordanRank}
Let $A_j=(E_{n_j}; x E_{n_j}+J_{n_j})$ be an $n_j\times n_j\times 2$ tensor
for $j=1,\ldots,\ell$
and $X$ an arbitrary $n'\times n'$ matrix.
Then
$$\rank_{\FFF}(\diag(A_1,\ldots,A_\ell,(E_{n'};X))) \geq 
  \sum_{j=1}^\ell n_j+n'+ \ell \enspace.$$
\end{thm}

\begin{proof}
It suffices to show the claim when $x=0$.
We take 
$$\diag(A_1,\ldots,A_\ell,(E_{n'};X))$$
 as an array with $n$ slices of $2 \times n$ matrices $B_1, B_2, \ldots, B_n$,
where 
$n=\sum_{j=1}^\ell n_j+n'$.
Since $\langle B_1,\ldots,B_n\rangle=n$, by applying 
Theorem~\ref{thm:x2mn>} for
$$A'_1=B_1, A'_2=B_{n_1+1},\ldots, A'_\ell=B_{n_1+\cdots+n_{\ell-1}+1}\enspace,$$
we can show the claim straightfowardly.
\end{proof}

\section{Decomposition and Rank\label{sec:companion}}

Now we recall that the maximal rank of tensors
with 2 slices was given by the following theorem.

\begin{thm}[{cf. \cite[Theorem~3.5]{JaJa:1979b}}]
\label{thm:thm3.5jaja79b}
$$\maxrank_{\FFF}(m,n,2)=
\min\left(n+\left\lfloor\frac{m}{2}\right\rfloor,m+\left\lfloor\frac{n}{2}\right\rfloor,2m,2n\right) \enspace.$$ 
\end{thm}
In this section 
we determine all tensors which attain the maximal rank.

First we consider about the rank of $(E_n;A)$.
\jaja{} discussed ranks by using invariant
polynomials \cite{JaJa:1979a,JaJa:1979b}.

Let $\KKK$ be an arbitrary field and $x$ an indeterminate over $\KKK$.
For a matrix $A(x)$ with entries in $\KKK[x]$, we denote by
$e_i(A(x))$ the $i$-th elementary divisor of $A(x)$.
If we denote the greatest common divisor of $i$-minors of $A(x)$ by
$d_i(A(x))$, then $e_i(A(x))=d_i(A(x))/d_{i-1}(A(x))$
in case $d_{i-1}(A(x))\neq 0$.

Here we recall a basic fact.

\begin{lemma}\label{lem:sim cri}
Let $A$, $B$ be $n\times n$ matrices with entries in $\KKK$.
Then $B$ is similar to $A$ if and only if
$$
e_i(xE_n-A)=e_i(xE_n-B)\qquad\text{for $i=1,2,\ldots,n$}\enspace.
$$
\end{lemma}
Note 
$e_1(xE_n-A)e_2(xE_n-A)\cdots e_n(xE_n-A)=\det(xE_n-A)\neq 0$
for an $n\times n$ matrix $A$ with entries in $\KKK$.
In particular,
$e_n(xE_n-A)\neq 0$.

Now we recall the result of \jaja. 
Let $A$ be an $n\times n$ matrix.
\jaja{} called $e_{n-i+1}(xE_n-A)$ the $i$-th invariant polynomial
of $A$ and denoted as $p_i(A)$.
\begin{thm}[{\cite[Theorem 3.3 and proof of Theorem 3.1]{JaJa:1979b}}]
\label{thm:jaja3.3}
Let $A$ be an $n\times n$ matrix and
$k$ the number of those $p_i(A)$'s which cannot be factored into 
distinct linear factors over $\KKK$.
Suppose $\card(\KKK)\geq \deg p_1(A)$.
Then
$\rank_\KKK(E_n;A)\leq n+k$.
In fact, $(E_n;A)$ is diagonalizable after adding 
$k$ tensors of rank $1$.
\end{thm}

The following example shows that Theorem~\ref{thm:thm3.5jaja79b} does not
hold over the Galois field $GF(2)$ and 
thus the condition $\card(\KKK)\geq \deg p_1(A)$
can not be removed in Theorem~\ref{thm:jaja3.3}.

\begin{prop}
For $A=\begin{pmatrix}
0&0&1\\
1&0&1\\
0&1&0
\end{pmatrix}$, it holds that
$
\rank_{GF(2)}(E_3;A)\geq 5 \enspace.
$
\end{prop}

\begin{proof}
Supposing that $\rank_{GF(2)}(E_3;A)\leq 4$ we show a contrary.
There are 
$\aaa_i$, $\bbb_i \in GF(2)^3$ and $\alpha_i$, $\beta_i \in GF(2)$
for $1\leq i\leq 4$ such that
$$
E_3=\sum_{i=1}^4\aaa_i\alpha_i\bbb_i^T,\quad
A=\sum_{i=1}^4\aaa_i\beta_i\bbb_i^T \enspace.
$$
Changing the suffix if necessary, we may assume that
$\aaa_1$, $\aaa_2$, $\aaa_3$ are linearly independent
and $\alpha_1$, $\alpha_2$, $\alpha_3\neq0$.
Since we are working over $GF(2)$,
this means $\alpha_1=\alpha_2=\alpha_3=1$.
On the other hand, since $(E_3;A)$ is not diagonalizable,
we see that $\bbb_4\neq \zerovec$.
And, since $\rank A=3$, by changing the suffix if necessary,
we may assume that
$\bbb_2$, $\bbb_3$, $\bbb_4$ are linearly independent
and $\beta_2$, $\beta_3$, $\beta_4\neq0$.
Again this implies that $\beta_2=\beta_3=\beta_4=1$.

Therefore, we see that
$$
E_3+A=E_3-A\in\langle\aaa_1\bbb_1^T, \aaa_4\bbb_4^T\rangle \enspace.
$$
This contradicts to the fact that $\rank(E_3+A)=3$.
\end{proof}

Note that
since $p_n(A)\mid p_{n-1}(A)\mid\cdots\mid p_2(A)\mid p_1(A)$,
$p_j(A)$ can be factored into distinct linear factors over $\KKK$
if and only if $j>k$, in the notation of Theorem \ref{thm:jaja3.3}.

\jaja{} \cite[Theorem 3.6]{JaJa:1979b} showed the reverse inequality on the
assumption that $p_1(A)$ can be factored into 
(not necessarily distinct) linear factors over $\KKK$.
Here we show the reverse inequality without any assumption.

\begin{thm}\label{thm:rev of jaja3.3}
Let $A$ and $k$ be as in Theorem \ref{thm:jaja3.3}.
Then
$$
\rank_\KKK(E_n;A)\geq n+k \enspace.
$$
\end{thm}

\begin{proof}
Set $\rank_\KKK(E_n;A)=n+q$.
We want to show that $q\geq k$, and 
so we may assume that $q<n$.
Take $\aaa_1, \ldots, \aaa_{n+q}$, $\bbb_1, \ldots, \bbb_{n+q}\in\KKK^n$
and $\alpha_1, \ldots, \alpha_{n+q}$, 
$\beta_1, \ldots, \beta_{n+q} \in \KKK$ 
such that
$$
E_n=\sum_{j=1}^{n+q}\aaa_i\alpha_i\bbb_i^T,\quad
A=\sum_{j=1}^{n+q}\aaa_i\beta_i\bbb_i^T \enspace.
$$
Changing the suffix if necessary, we may assume that
$\aaa_1, \ldots, \aaa_n$ are linearly independent and
$\alpha_1, \ldots, \alpha_n \ne 0$ since $\rank(E_n)=n$.
By exchanging $\alpha_i\aaa_i$ by $\aaa_i$ for $1\leq i\leq n$, 
we may assume $\alpha_1=\cdots=\alpha_n=1$.
Set $\dim\langle\bbb_{n+1},\ldots,\bbb_{n+q}\rangle=q'$.
Then by changing the suffix within 
$\{n+1, \ldots, n+q\}$
if necessary,
we may assume that
$\bbb_{n+1}$, \ldots,  $\bbb_{n+q'}$ is a basis of 
$\langle\bbb_{n+1},\ldots,\bbb_{n+q}\rangle$.
Then
$\bbb_j\in\langle\bbb_{n+1},\ldots,\bbb_{n+q'}\rangle$
for $j>n+q'$.
Since $\dim\langle\bbb_1$, \ldots, $\bbb_n$, \ldots, $\bbb_{n+q}\rangle=n$,
we may further assume, 
by changing the suffix 
within $\{1, \ldots, n\}$
if necessary, that 
$\bbb_{q'+1}$, \ldots, $\bbb_n$, \ldots,  $\bbb_{n+q'}$ are 
linearly independent.

Then there are nonsingular matrices $P$ and $Q$ with entries in $\KKK$
such that
$$
P(\aaa_1,\ldots,\aaa_{n+q})=(E_n,\ast),\qquad
\begin{pmatrix}
\bbb_1^T\\
\bbb_2^T\\
\vdots\\
\bbb_{n+q}^T
\end{pmatrix}Q
=
\begin{pmatrix}
\ast&\ast\\
E_{n-q'}&O\\
O&E_{q'}\\
O&\ast
\end{pmatrix} \enspace.
$$

Since 
\begin{equation*}
\begin{split}
xE_n-A & =\sum_{j=1}^{n+q}\aaa_i(\alpha_i x- \beta_i)\bbb_i^T \\
 & =(\aaa_1,\ldots,\aaa_{n+q})
    \diag(\alpha_1 x-\beta_1,\ldots, \alpha_{n+q}x-\beta_{n+q})
\begin{pmatrix}
\bbb_1^T\\
\bbb_2^T\\
\vdots\\
\bbb_{n+q}^T
\end{pmatrix}\enspace,
\end{split}
\end{equation*}
we wee that
\begin{eqnarray*}
P(xE_n-A)Q&=&
(E_n,\ast)\diag(\alpha_{1}x-\beta_1,\ldots,\alpha_{n+q}x-\beta_{n+q})
\begin{pmatrix}
\ast&\ast\\
E_{n-q'}&O\\
O&E_{q'}\\
O&\ast
\end{pmatrix}\\
&=&
\begin{pmatrix}
\ast&\ast\\
\diag(x-\beta_{q'+1},\cdots,x-\beta_n) &\ast
\end{pmatrix} \enspace.
\end{eqnarray*}
Therefore
$d_{n-q'}(xE_n-A)=d_{n-q'}(P(xE_n-A)Q)$
divides $\prod_{j=q'+1}^n(x-\beta_j)$ and can be factored into linear 
factors over $\KKK$.
Since $q\geq q'$ and
$p_{q+1}(A)=e_{n-q}(xE_n-A)$ divides $d_{n-q}(xE-A)$,
we see that $p_{q+1}(A)$ can be factored into linear factors over $\KKK$.

By assumption, $p_k(A)$ cannot be factored into distinct linear factors 
over $\KKK$.
So $p_k(A)$ has an irreducible factor of degree greater than 1 and/or
$p_k(A)$ has a multiple linear factor.

In the first case, $q+1>k$ since $p_{q+1}(A)$ does not have an 
irreducible factor whose degree is greater than $1$.
Therefore $q\geq k$.
Now assume that $(x-\beta)^2$ divides $p_k(A)$. 
Then $A$ is similar to
$B=\diag(\beta E_{m_1}+J_{m_1},\ldots, \beta E_{m_{k}}+J_{m_k},A')$
for appropriate $A'$ by Lemma \ref{lem:sim cri}.
Therefore
$$
\rank_\KKK(E_n;A)=\rank_\KKK(E_n;B)\geq n+k
$$
by Theorem \ref{thm:JordanRank}.
\end{proof}

As a corollary, we obtain the main theorem of this section.

\begin{thm}\label{thm:companion}
Let $A$ be an $n\times n$ matrix and let $\alpha_{\FFF}(A,x)$
be the number of 
Jordan blocks whose sizes are greater than or equal to
$2$ for an eigenvalue $x$ of $A$.  
Then 
$$
\rank_{\FFF}(E_n;A)= n + \max_x \alpha_{\FFF}(A,x) \enspace,
$$
where we treat $C_k(c,s)+J_k\otimes E_2$ as a Jordan block of size $2k$
if $\FFF=\RRR$.
Furthermore, the tensor $(E_n;A)$ is diagonalizable after adding
$\max_x \alpha_{\FFF}(A,x)$ tensors of rank $1$.
\end{thm}


Then we have easily to obtain a border rank.

\begin{prop}[{\cite[Proposition~3.3]{Bini:1980}}]
For a border rank $\brk_{\FFF}(E_n;A)$, we have
$$\brk_{\CCC}(E_n;A)=n \text{ and } \brk_{\RRR}(E_n;A) =n,n+1\enspace.$$
In particular,
$$\maxbrk_{\CCC}(n,n,2)=n \text{ and } \maxbrk_{\RRR}(n,n,2)=n+1\enspace.$$
\end{prop}

\begin{proof}
There is a sequence $\{A_j\}$ of $n\times n$ matrices whose
eigenvalues in $\CCC$ are distinct each other and converges to $A$.
Then $\rank_{\CCC}(E_n,A_j)=n$ and $\rank_{\RRR}(E_n,A_j)=n,n+1$
for each $j$.
If $A$ has a complex, not real eigenvalues, then $A_j$ has also
for sufficiently large $j$ and thus $\rank_{\RRR}(E_n,A_j)=n+1$.
\par
For arbitrary $n\times n\times 2$ tensor $(X;Y)$, there is a sequence
$\{(X_j;Y_j)\}$ such that $X_j$ is nonsingular and eigenvalues
of $X_j^{-1}Y_j$ are distinct each other for each $j$.
Thus the claim follows from 
$\rank_{\FFF}(X_j;Y_j)=\rank_{\FFF}(E_n;X_j^{-1}Y_j)$.
\end{proof}

Two tensors $T$ and $T'$ are called {\sl equivalent} if there are nonsingular matrices $P$ and $Q$ such that $PTQ=T'$.

Before closing this section we show the rank of
a tensor $(A;B)$ having  a Kronecker canonical form.
Let $A$ and $B$ be $m\times n$ rectangular matrices.
The rank of a tensor $(A;B)$ is obtained by 
its Kronecker canonical form
(cf. \cite[Theorem~5]{JaJa:1979b}).
If $(A;B)$ is equivalent to one consisting of the direct sum of 
an $m_A\times n_A\times 2$ tensor $(O;O)$ of type \eqref{typeA},
an $m_E^{(i)}\times (m_E^{(i)}+1)\times 2$ tensor of type \eqref{typeE} for $1\leq i\leq \ell_E$, 
and an $(n_F^{(i)}+1)\times n_F^{(i)}\times 2$ tensor of type \eqref{typeF} for $1\leq i\leq \ell_F$,
and tensors of type \eqref{typeB}, \eqref{typeD} and
in addition if $\FFF=\RRR$, tensors of type \eqref{typeC}.
Let $\alpha$ be the maximal integer among
the number of $(x E_k+J_k;E_k)$ of type \eqref{typeB} with $k\geq 2$ for each $x$,
the number of $(E_k;J_k)$ of type \eqref{typeD} with $k\geq 2$,
and in addition if $\FFF=\RRR$
the number of $(C_k(c,s)+J_k\otimes E_j;E_{2k})$ with $k\geq 1$ for each $(c,s)$, $s\ne 0$.
Put $m_E=\sum_{i=1}^{\ell_E} m_E^{(i)}$ and
$n_F=\sum_{i=1}^{\ell_F} n_F^{(i)}$ for short.

\begin{thm}  \label{thm:diagonal}
It holds 
$m-m_A+\ell_E=n-n_A+\ell_F$
and
$$\rank_{\FFF}(A;B) = \alpha + m-m_A+\ell_E\enspace.$$
In fact there is a tensor $T$ of rank $\alpha+\ell_E+\ell_F$ such that
$(A;B)+T$ is diagonalizable.
\end{thm}

\begin{proof}
We may assume that $(A;B)$ is of a Kronecker canonical form.
Let 
$$
(A;B)=\diag((O;O),(A_1;B_1),(A_2;B_2)),
$$ 
where
$(A_1;B_1)$ is an $(m_E+n_F+\ell_F)\times(m_E+n_F+\ell_E)\times 2$ tensor
consisting of tensors of type \eqref{typeE} and \eqref{typeF} and
$(A_2;B_2)$ is a tensor consisting of tensors of type \eqref{typeB}, 
\eqref{typeC} and \eqref{typeD}.
By Lemma~\ref{lem:GSDcaseE&F} the tensor $(A_1;B_1)$ is diagonalizable
after adding at most $\ell_E+\ell_F$ tensors of rank $1$.
Since a tensor of type \eqref{typeB}, \eqref{typeC}, \eqref{typeD} consists $2$ slices of square matrices, we have $(m-m_A)-(n-n_A)=\ell_F-\ell_E$.
For simplicity, let 
$$p=m-m_A-m_E-n_F-\ell_F$$ 
which is the size of the square matrix $A_2$.
Take $d \in \FFF$ so that $A_2+d B_2$ is nonsingular.
Direct summands of $(A_2;B_2)$ are 
$1$ to $1$ corresponding to Jordan blocks of $(A_2+d B_2)^{-1}B_2$.
Furthermore, Jordan blocks with eigenvalue $0$ come from 
tensors of type \eqref{typeD}, and 
if $\FFF=\RRR$ Jordan blocks with non-real eigenvalues come from tensors of type \eqref{typeC}.
Thus $\alpha=\max_x \alpha_{\FFF}((A_2+dB_2)^{-1}B_2,x)$.
By Theorem~\ref{thm:companion} $(A_2+dB_2;B_2)$ and then $(A_2;B_2)$
is diagonalizable after adding $\alpha$ tensors of rank $1$.
Therefore $(A;B)$ is diagonalizable after adding a tensor
of rank at most $\alpha+\ell_E+\ell_F$ and the rank of the obtained diagonal
tensor is equal to $p+m_E+n_F=m-m_A-\ell_F$.
Moreover, it follows by Theorem~\ref{thm:companion} and \cite[Theorem~2.4]{JaJa:1979b} 
that
\begin{equation*}
\begin{split}
\rank_{\FFF}(A;B) & = \rank_{\FFF}(A_2;B_2)
  +(m_E+\ell_E)+(n_F+\ell_F) \\
 & = \rank_{\FFF}(E_p;(A_2+dB_2)^{-1}B_2) +m-m_A-p+\ell_E  \\
 & = \alpha + m-m_A + \ell_E \enspace. \\
\end{split}
\end{equation*}
\end{proof}

As a corollary, we obtain all Kronecker canonical forms 
giving the maximal rank.
We denote by $X^{\oplus k}$ the direct sum of $k$ copies of a tensor $X$.

\begin{cor}
Suppose $m\leq n\leq 2m$ and $\rank_{\FFF}(A;B)=\maxrank_{\FFF}(m,n,2)$.
If $n$ is even, then $(A;B)$ is equivalent to 
$$\diag(Y^{\oplus\alpha},((0,1);(1,0))^{\oplus \ell_E})$$
and otherwise $(A;B)$ is equivalent to one of the following tensors:
\begin{itemize}
\item[{\rm (i)}] $\diag(Y^{\oplus\alpha},((0,1);(1,0))^{\oplus \ell_E}),\zerovec)$
\item[{\rm (ii)}] $\diag(Y^{\oplus\alpha},((0,1);(1,0))^{\oplus \ell_E},((0,1);(1,0))^T)$
\item[{\rm (iii)}] $\diag(Y^{\oplus\alpha},((0,1);(1,0))^{\oplus \ell_E},(x;1))$
\item[{\rm (iv)}] $\diag(Y^{\oplus\alpha},((0,1);(1,0))^{\oplus \ell_E},(1;0))$
\item[{\rm (v)}] $\diag((xE_2+J_2;E_2)^{\oplus(\alpha-1)},(xE_3+J_3;E_3),((0,1);(1,0))^{\oplus \ell_E})$ 
\item[{\rm (vi)}]
$\diag((E_2;J_2)^{\oplus(\alpha-1)},(E_3;J_3),((0,1);(1,0))^{\oplus \ell_E})$
\item[{\rm (vii)}] 
$\diag(Y^{\oplus\alpha},((0,1);(1,0))^{\oplus (\ell_E-1)},
  ((\zerovec,E_2);(E_2,\zerovec)))$
\end{itemize}
where $Y$ is $(xE_2+J_2;E_2)$, $(E_2;J_2)$, or $(C_1(c,s);E_2)$.
\end{cor}

\begin{proof}
As in the proof of Theorem~\ref{thm:diagonal}, put 
$p=m-m_A-m_E-n_F-\ell_F$ and 
let $(A_2;B_2)$ be the tensor consisting of all direct summands 
of type \eqref{typeB}, \eqref{typeC} and \eqref{typeD} in $(A;B)$.
Then $m=m_A+p+m_E+n_F+\ell_F$ and $n=n_A+p+m_E+n_F+\ell_E$.
Since 
$\lfloor n/2\rfloor=\rank_{\FFF}(A;B)-m$, 
it holds
$$(p-2\alpha)+(m_E-\ell_E)+2m_A+n_A+n_F=n-2\lfloor n/2\rfloor\enspace.$$
Note that $p\geq 2\alpha$, $m_E\geq \ell_E$ and $n_F\geq \ell_F$.
Then $m_A=0$.
If $n$ is even it holds that $p=2\alpha$, $m_E=\ell_E$, $n_A=n_F=0$.
$p=2\alpha$ yields that 
$(A_2;B_2)=Y^{\oplus\alpha}$ for some 
$Y=(xE_2+J_2;E_2), (E_2;J_2), (C_1(c,s);E_2)$ and
$m_E=\ell_E$ implies that the direct summand of tensors of type \eqref{typeE}
is $((0,1);(1,0))^{\oplus \ell_E}$.
Therefore $(A;B)$ is equivalent to 
$\diag(Y^{\oplus\alpha},((0,1);(1,0))^{\oplus \ell_E})$
when $n$ is even.
Now let $n$ be odd.
One of $p-2\alpha$, $m_E-\ell_E$, $n_A$ and $n_F$ is one and the others
are all zero.
The tensor $(A;B)$ is equivalent to the tensor \eqref{enu:1}
if $n_A=1$ and to the tensor \eqref{enu:2} if $n_F=1$.
In the case when $p=2\alpha+1$, $(A;B)$ is equivalent to 
\eqref{enu:3}, \eqref{enu:4}, \eqref{enu:4-a} or \eqref{enu:4-b}.
Finally if $m_E=\ell_E+1$, then $(A;B)$ is equivalent to \eqref{enu:5}.
\end{proof}

\begin{cor}
Let $m$ and $n$ be positive integers with $m\leq n$.
Any $m\times n\times 2$ tensor is diagonalizable after adding
at most $\lfloor n/2\rfloor$ tensors of rank $1$.
\end{cor}

\begin{proof}
Let $L_k=((\zerovec,E_k);(E_k,\zerovec))$ be a $k\times (k+1)\times 2$
tensor of type \eqref{typeE}. Then $L_k^T$ is a $(k+1)\times k\times 2$
tensor of type \eqref{typeF}.
By Lemma~\ref{lem:GSDcaseE&F}, for a tensor $\diag(L_a,L_b^T)$,
$\diag(L_a,L_b^T)+T$ is diagonalizable for some tensor $T$ of rank
$2$.  In particular, if $a,b>0$ and $a+b\geq 3$, then 
$\diag(L_a,L_b^T)$ is diagonalizable after adding some tensor of rank
at most $\lfloor(a+b+1)/2\rfloor$.
We show $\diag(L_1,L_1^T)$ is diagonalizable after adding some tensor of rank
$1$.
Set $\diag(L_1,L_1^T)=(X;Y)$ 
and $M=\begin{pmatrix} 0&0&0\cr 1&1&0\cr 0&0&0\end{pmatrix}$.
Then $X+M$ is nonsingular and $(X+M)^{-1}(Y+M)$ has eigenvalues $\pm1,0$.
Thus $(X;Y)+(M;M)$ is diagonalizable.
Therefore for $a,b>0$, the $(a+b+1)\times (a+b+1)\times 2$ tensor
$\diag(L_a,L_b^T)$ is diagonalizable 
after adding adequate tensor of rank at most $\lfloor(a+b+1)/2\rfloor$.
\par
Suppose that $\ell_E\geq \ell_F$.
Then by Lemma~\ref{lem:BCD} and the above observation,
we see that $(A;B)$ is a direct sum of tensors 
each one is diagonalizable after adding a rank 1 tensor and 
has at least 2 columns.
So the result follows.
We can treat the case where $\ell_F\geq\ell_E$ by the same way.
So
we complete the proof.
\end{proof}

\end{document}